\documentclass[11pt]{article}

\def\d{\delta}

\def\s{\sigma}

\def\ot{\otimes}

\def\lra{\longrightarrow}

\def\mapright#1{\smash{\mathop{\lra}\limits^{#1}}}

\def\NM{{\mathbf N}}

\newtheorem{thm}{Theorem}

\newtheorem{prop}{Proposition}

\begin{document}

\title{\bf Brace algebras and the cohomology comparison theorem$^{(*)}$.}

\author{\bf Fr\'ed\'eric Patras$^{(**)}$}

\date{ }
\maketitle

\vspace{ 0,6 cm}

       \noindent {\bf Abstract.} {\small The Gerstenhaber and Schack
       cohomology comparison
theorem asserts that there is a cochain equivalence between
the Hochschild complex of a certain algebra and the usual singular cochain complex of a space. We show that this comparison theorem preserves the brace algebra structures.
This result gives a structural reason for the recent results
establishing fine topological structures on the Hochschild cohomology,
and a simple way to derive them from the corresponding properties of
cochain complexes.

\vspace{ 0,4 cm}

\noindent{\bf A.M.S Classification.} 16E40; 55N10; 18D50; 55P48

\vspace{ 0,4 cm}


{\vspace{3,2 cm}

\rule{12 cm}{0,2mm}
\vspace{ 0,2 cm}
{\footnotesize{
(*) Revised version of : The bar construction as a Hopf algebra,
December 2001.\\
(**)  CNRS UMR 6621 - Universit\'e de Nice,
Math\'ematiques, Parc Valrose,
06108 Nice cedex 2,
France.  {\it patras@math.unice.fr} \ Tel. 00 33 (0) 492076262, Fax. 00 33 (0) 493517974.}}
\eject
\section*{Introduction}
A theorem of Gerstenhaber and Schack (the cohomology comparison
theorem, CCT) asserts that, for a given triangulated topological space,
there exists an associative algebra $A$ and a quasi-isomorphism of
cochain complexes between the cohomological Hochschild complex of $A$
and the singular cochain complex of the space \cite{gs1,gs2}. Besides,
there are brace differential graded algebra (BDGA) structures on the
cohomological Hochschild complexes of associative algebras and on the
singular cochain complexes. We prove that the Gerstenhaber-Schack
quasi-isomorphism preserves these algebraic structures. This result
should make clear the origin of the fine topological structures
appearing on the Hochschild cohomology according to the Deligne
conjecture \cite{ks,t,ms}, and why cochain algebras and Hochschild
complexes share many algebraic properties, the CCT providing a
systematic tool for "structure transportation" between the two
theories.

\section{Brace differential graded algebras}
Let us introduce first BDGAs. These algebras first appeared in the
work of Getzler-Jones on algebras up to homotopy (without a
specific name) as a particular case of $B_\infty$-algebras, associated
in particular to Hochschild complexes of associative algebras, see
\cite[Sect. 5.2]{gj}. When Gerstenhaber and Voronov studied them more
in detail \cite{gv,vg,v}, they decided to call these algebras homotopy
$G$-algebras. However, this
terminology appeared to be a misleading one after Tamarkin had shown
that the name $G$(erstenhaber)-algebra up to homotopy should be
naturally given to another class of algebras \cite{t}. We therefore
call them by a name that reflects their properties and
should not create confusion, namely: brace differential graded
algebras.\par The basic idea is that BDGAs are
associative differential graded algebras together with extra (brace)
operations that behave exactly as the Kadeishvili-Getzler brace
operations on the Hochschild cohomological complex of an associative
algebra \cite{k,ge}. We write, as usual, $B(A)$ for the cobar coalgebra
over a differential graded algebra (DGA) $A$, where the product is
written $\cdot$  and the differential (of degree +1) $\d$. That is,
$B(A)$ is the cofree graded coalgebra
$T(A[1]):=\bigoplus\limits_{n\in\NM}A[1]^{\ot n}$ over the desuspension
$A[1]$ of $A$ ($A[1]_n:=A_{n+1}$). We use the bar notation and write
$[a_1|...|a_n]$ for $a_1\ot ...\ot a_n\in A[1]^{\ot n}$. In particular,
the coproduct on $T(A[1])$ is given by: $$\Delta
[a_1|...|a_n]:=\sum\limits_{i=0}^n[a_1|...|a_i]\ot [a_{i+1}|...|a_n].$$
There is a differential coalgebra structure on $B(A)$ induced by the
DGA structure on $A$. In fact, since $B(A)$ is cofree as a graded
coalgebra, the properties of the cofree coalgebra functor imply that,
in general, a coderivation $D\in Coder(B(A))$ is entirely determined by
the composition (written as a degree 0 morphism): $$\tilde D:\
B(A)\mapright{D}B(A)[1]\mapright{p}A[2],$$ where $p$ is the natural
projection. In particular, the differential $d$ on $B(A)$ is induced by
the maps: $$\d :A[1]\mapright{} A[2],$$ and $$\mu :A[1]\ot
A[1]\mapright{} A[2],$$ where $\mu (a,b):=(-1)^{|a|}a\cdot b$. The
algebra $A$ is a BDGA if it is provided with a set of
extra-operations called the braces: $$B_k: A[1]\ot A[1]^{\ot
k}\mapright{} A[1],\ k\geq 1,$$ satisfying certain relations. These
relations express exactly the fact that the braces have to induce a
differential Hopf algebra structure on $B(A)$. Explicitly, the
relations satisfied by the braces are then \cite[Sect. 5.2]{gj} and
\cite{kh,v} (we use Getzler's notation:
$v\{v_1,...,v_n\}:=B_n(v\ot (v_1\ot
...\ot v_n))$)  :
\begin{enumerate}
\item The brace relations
(the associativity relations for the product on $B(A)$).
$$(v\{v_1,...,v_m\})\{w_1,...,w_n\}=\sum\limits_{0\leq i_1\leq
j_1\leq ...\leq i_m\leq j_m\leq
n}(-1)^{\sum\limits_{k=1}^m(|v_k|-1)(\sum\limits_{l=1}^{i_k}(|w_l|-1))}$$
$$v\{w_1,...,w_{i_1},v_1\{w_{i_1+1},...,w_{j_1}\},w_{j_1+1},...,
v_m\{w_{i_m+1},...,w_{j_m}\},w_{j_m+1},...,w_n\},$$ with the usual
conventions on indices: for example, an expression such as
$v_5\{w_7,...,w_6\}$ has to be read $v_5\{\emptyset\}=v_5$. \item The
distributivity relations of the product w.r. to the braces. $$(v\cdot
w)\{ v_1,...,v_n\}=\sum\limits_{k=0}^n(-1)^{|
w|\sum\limits_{p=1}^k(|v_p| -1)}v\{v_1,...,v_k\}\cdot
w\{v_{k+1},...,v_n\},$$
\item The boundary relations.
$$\d (v\{v_1,...,v_n\})-\d v\{v_1,...,v_n\}$$
$$+\sum\limits_{i=1}^n(-1)^{|v|+|v_1|+...+|v_{i-1}|-i+1}v\{v
_1,...,\d v_i,...,v_n\}$$ $$=(-1)^{|v|(|v_1|-1)}v_1\cdot (
v\{v_2,...,v_n\})$$
$$-\sum\limits_{i=1}^{n-1}(-1)^{|v|+|v_1|+... + | v_{ i
}|-i-1}v\{v_1,...,v_i\cdot v_{i+1},...,v_n\}$$
$$+(-1)^{|v|+|v_1|+...+| v_{ n-1 }|-n}(v\{v_1,...,v_{n-1}\})\cdot
v_n.$$

\end{enumerate}

Let us
write down explicit formulas for the BDGA structure on
the Hochschild cochain complex $C^{\ast}(A,A)$ of an associative
algebra $A$ over a commutative unital ring $k$. Recall that
$C^n(A,A)=Hom_k(A^{\ot n},A)$ and that the brace operations on
$C^{\ast}(A,A)$ are the multilinear operators defined for
$x,x_1,...,x_n$ homogeneous elements in $C^{\ast}(A,A)$ and
$a_1,...,a_m$ elements of $A$ by:\\
$$\{x\}\{x_1,...,x_n\}(a_1,...,a_m):=\sum\limits_{0\leq i_1\leq
i_1+|x_1|\leq i_2\leq ...\leq i_n+|x_n|\leq
n}(-1)^{\sum\limits_{k=1}^ni_k\cdot (|x_k|-1)}$$
$$x(a_1,...,a_{i_1},x_1(a_{i_1+1},...,a_{i_1+|x_1|}),...,a_{i_n},x_n(a_
{ i _n+1},...,a_{i_n+|x_n|}),...a_m).$$
The other operations defining the BDGA structure, $\d$
and $\cdot$ are, respectively, the
Hochschild coboundary and the cup product.\par
There is also
a BDGA structure on the
cochain complex of a simplicial set \cite{gj,gv}.\par
Recall that a simplicial set is a contravariant functor
from the category $\bf\Delta$ of finite sets $[n]=\{0,...,n\}$ and
increasing morphisms to $\bf Sets$. For a simplicial set $S:{\bf
\Delta}\lra \bf Sets$, for $\sigma\in S_n:=S([n])$, and for a strictly
increasing sequence $0\leq a_0< ...<a_m\leq n$, we write $\sigma
(a_0,...,a_m)$ for $i_a^\ast (\sigma )\in S_m$, where $i_a$ is the
unique map from $[m]$ to $[n]$ sending $[m]$ to $\{a_0,...,a_m\}$.
Define a map $\Delta_{1,r}$ from the singular complex of $X$, $C_\ast
(X)$ to $C_\ast (X)\ot C_\ast (X)^{\ot r}$ as follows. For $\sigma\in X_n$, set:\\
$$\Delta_{1,r}(\sigma ):=\sum\limits_{0\leq b_1'\leq b_1\leq ...\leq
b_r'\leq b_r\leq n}(-1)^{\sum\limits_{k=1}^r((b_k-b_k')b_k')}$$
$$\sigma
(0,1,...,b_1',b_1,b_1+1,...,b_2',b_2,...,b_r',b_r,...,n-1,n)$$ $$\ot
(\sigma (b_1',...,b_1)\ot \sigma (b_2',...,b_2)\ot ...\ot \sigma
(b_r',...,b_r)).$$
Dualizing $\Delta_{1,r}$, we get a
map from $C^\ast (X)\ot C^\ast (X)^{\ot r}$ to $C^\ast
(X)$. By analogy with the case of Hochschild cochains, we write
$\s\{\s_1,...,\s_r\}$ for $\Delta_{1,r}^\ast (\s\ot (\s_1\ot ...\ot
\s_r))$. These brace operations on cochains, together with the
simplicial coboundary and the cup product induce a BDGA
structure on the bar coalgebra on $C^\ast (X)$.

\section{The cohomology comparison theorem}
Recall Gerstenhaber and Schack's cohomology comparison theorem
for finite simplicial complexes and their incidence algebras
\cite{gs1,gs2}. Let $\Sigma$ be a subset of the set of subsets of a
finite set $S$. Then, $\Sigma$ is a finite simplicial complex if, for
each $\sigma$ in $\Sigma$, the set of subsets of $\sigma$ (viewed as a
subset of $S$) is a subset of $\Sigma$. If $|\sigma|=n+1$, $\sigma$ is
called a $n$-simplex of $\Sigma$. The elements of $\Sigma$ are ordered
by inclusion and, in particular, we can view $\Sigma$ as a poset. We
write $\leq$ for the inclusion between the
simplices of $\Sigma$.\par To $\Sigma$ are associated canonically two
objects, both of which compute its simplicial cohomology. The first
one, written $\hat\Sigma$, is the usual barycentric subdivision of
$\Sigma$. It is the simplicial set whose $n$-simplices are the ordered
morphisms (weakly increasing maps) from the ordered set $\{0,...,n\}$
to $\Sigma$ or, equivalently, the increasing sequences in $\Sigma$,
written $\s_0\leq ...\leq \s_n$. The cohomology of $\hat\Sigma$ is
isomorphic to the simplicial cohomology of $\Sigma$.\par The other
object is the incidence algebra $I_\Sigma$ of the poset $\Sigma$: it is
the algebra generated linearly (over a commutative ring $k$) by the
pairs of simplices $(\sigma ,\sigma ')$ with $\sigma\leq \sigma '$. The
product of two pairs $(\sigma ,\sigma ')$ and $(\beta ,\beta ')$ is
$(\sigma ,\beta ')$ if $\sigma '= \beta$ and $0$ else. This algebra is
a triangular algebra. The Hochschild cohomology of such algebras can be
computed explicitly by means of a spectral sequence, introduced
recently by S. Dourlens \cite{d}. We refer from now on to
\cite{gs2} and \cite{d} for the general properties of the Hochschild
cohomology of triangular and incidence algebras that are recalled
below.\par The incidence algebra $I_\Sigma$ has a separable subalgebra
$S_\Sigma$ generated as a $k$-algebra by the pairs $(\sigma ,\sigma )$.
The $n$-cochains of the Hochschild complex of $I_\Sigma$ relative to
$S_\Sigma$ are the elements of $Hom_{S_\Sigma
-S_\Sigma}(I_\Sigma^{\ot_{S_\Sigma}n},I_\Sigma )$, with the usual
formula for the Hochschild coboundary. This relative Hochschild
complex, written $C_{S_\Sigma}^\ast (I_\Sigma ,I_\Sigma )$, computes
the usual Hochschild cohomology of $I_\Sigma$. A direct inspection
shows that $Hom_{S_\Sigma
-S_\Sigma}(I_\Sigma^{\ot_{S_\Sigma}n},I_\Sigma )$ is generated linearly
(over $k$) by the maps sending a given tensor product
$((\sigma_0,\sigma_1),(\sigma_1,\sigma_2),...,(\sigma_{n-1},\sigma_n))$
to $(\sigma_0,\sigma_n)$, where the $(\sigma_i,\sigma_{i+1})$s belong
to the set of generators of $I_\Sigma$, and all the other tensor
products of generators of $I_\Sigma$ to $0$. \par The Gerstenhaber and
Schack cohomology comparison theorem states that there is a canonical
cochain isomorphism between the singular cohomology of the barycentric
subdivision of $\Sigma$ and this relative Hochschild complex. See
\cite{gs2} (in particular sect. 23) for details, generalizations, and a
survey of the history of this theorem. \begin{thm} There is a cochain
algebra isomorphism $\iota$ between the singular complex of
$\hat\Sigma$ and the relative Hochschild complex of $I_\Sigma$ given
by: for $f\in Hom_k(\hat\Sigma_n ,k)$ $$\iota
(f)((\sigma_0,\sigma_1),(\sigma_1,\sigma_2),...,(\sigma_{n-1},\sigma_n)
):=f(\sigma_0\leq \sigma_1\leq ...\leq \sigma_n)\cdot
(\sigma_0,\sigma_n).$$ In particular: $$HH^\ast (I_\Sigma ,I_\Sigma
)\cong H^\ast (\Sigma ,k).$$ \end{thm}

\begin{prop}The isomorphism
$\iota$ commutes with the action of the brace operations on $C^\ast
(\hat\Sigma ,k)$ and $C_{S_\Sigma}^\ast (I_\Sigma ,I_\Sigma
)$.\end{prop}
Indeed, let $f,f_1,...,f_k$ belong respectively to
$C^n(\hat\Sigma )$, $C^{n_1}(\hat\Sigma )$,..., $C^{n_k}(\hat\Sigma )$.
Let
$(\sigma_0\leq
\sigma_1\leq \sigma_2\leq ...\leq \sigma_{m-1}\leq \sigma_m)\in
\hat\Sigma_m$, where $m:=n+n_1+...+n_k-k $. Let
us also introduce the following useful convention. Let e.g.
$(\sigma_{i_0},\sigma_{i_1},k_1,k_2,\sigma_{i_3},...,\sigma_{i_q},k_p)$
be any sequence, the elements of which are either scalars, either
simplices of $\Sigma$, and assume that
$(\sigma_{i_0}\leq\sigma_{i_1}\leq ...\leq\sigma_{i_q})$ is a simplex
of $\hat\Sigma$. Then, we write
$f(\sigma_{i_0},\sigma_{i_1},k_1,k_2,\sigma_{i_3},...,\sigma_{i_q},k_p)
$ for $(\prod_{i=1}^pk_i)\cdot f(\sigma_{i_0}\leq\sigma_{i_1}\leq
...\leq\sigma_{i_q})$.\par Then, we have, according to the definition
of the braces:\par $f\{f_1,...,f_k\}(\sigma_0\leq\sigma_1\leq
...\leq\sigma_m)$ $$=\sum\pm
f(\sigma_0,...,\sigma_{i_1},f_1(\sigma_{i_1}\leq
...\leq\s_{i_1+n_1}),\s_{i_1+n_1},...,\s_{i_k},$$
$$f_k(\sigma_{i_k}\leq ...\leq\s_{i_k+n_k}),\s_{i_k+n_k},...,\s_m).$$
Therefore:\par $\iota
(f\{f_1,...,f_k\})((\sigma_0,\sigma_1),(\sigma_1,\sigma_2),...,(\sigma_
{m-1},\s_m))$ $$=\{\sum\pm
f(\sigma_0,...,\sigma_{i_1},f_1(\sigma_{i_1}\leq
...\leq\s_{i_1+n_1}),\s_{i_1+n_1},...,\s_{i_k},$$
$$f_k(\sigma_{i_k}\leq
...\leq\s_{i_k+n_k}),\s_{i_k+n_k},...,\s_m)\}\cdot (\s_0,\s_m).$$
$$=\sum\pm \iota
(f)((\sigma_0,\sigma_1),...,(\sigma_{i_1-1},\sigma_{i_1}),f_1(\sigma_{i
_1}\leq ...\leq\s_{i_1+n_1})\cdot (\sigma_{i_1},\sigma_{i_1+n_1}),$$
$$(\sigma_{i_1+n_1},\sigma_{i_1+n_1+1}),...
,(\sigma_{i_k-1},\sigma_{i_k}),$$ $$f_k(\sigma_{i_k}\leq
...\leq\s_{i_k+n_k})\cdot
(\sigma_{i_k},\sigma_{i_k+n_k}),...,(\s_{m-1},\s_m))$$ $$=\sum\pm \iota
(f)\{\iota (f_1),...,\iota
(f_k)\}((\sigma_0,\sigma_1),(\sigma_1,\sigma_2),...,(\sigma_{m-1},\s_m)
),$$ and the proof of the proposition follows.\par Notice that the
signs in the definition of the brace operations on cochains have been
chosen in such a way that the last identity holds.

\begin{thm}The morphism $\iota$ is an isomorphism of
BDGAs between the singular
cochain complex of the barycentric subdivision of a finite simplicial
complex $\Sigma$ and the $S_\Sigma$-relative Hochschild cochain complex
of the incidence algebra of $\Sigma$.\par In particular, the cohomology
comparison theorem of Gerstenhaber and Schack relating singular
cohomology and Hochschild cohomology can be realized, at the cochain
level, as a quasi-isomorphism of BDGAs.\end{thm}

 \end{document}